\documentclass[reqno,12pt]{amsart}

\usepackage{amscd,amssymb}

\theoremstyle{plain}
\newtheorem{Thm}[subsection]{Theorem}
\newtheorem{Cor}[subsection]{Corollary}
\newtheorem{Lem}[subsection]{Lemma}
\newtheorem{Prop}[subsection]{Proposition}

\theoremstyle{definition}
\newtheorem{Def}[subsection]{Definition}

\theoremstyle{remark}

\newtheorem{Rem}[subsection]{Remark}

\numberwithin{equation}{section}
\renewcommand{\rm}{\normalshape}

\newif\ifShowLabels
\ShowLabelstrue
\newdimen\theight
\def\TeXref#1{%
	\leavevmode\vadjust{\setbox0=\hbox{{\tt
		\quad\quad  {\small \rm #1}}}%
	\theight=\ht0
	\advance\theight by \lineskip
	\kern -\theight \vbox to
	\theight{\rightline{\rlap{\box0}}%
	\vss}%
	}}%

\ShowLabelsfalse

\renewcommand{\sec}[2]{\section{#2}\label{S:#1}%
	\ifShowLabels \TeXref{{S:#1}} \fi}
\newcommand{\ssec}[2]{\subsection{#2}\label{SS:#1}%
	\ifShowLabels \TeXref{{SS:#1}} \fi}

\newcommand{\sssec}[2]{\subsubsection{#2}\label{SSS:#1}%
	\ifShowLabels \TeXref{{SSS:#1}} \fi}

\newcommand{\refs}[1]{Section ~\ref{S:#1}}
\newcommand{\refss}[1]{Section ~\ref{SS:#1}}

\newcommand{\reft}[1]{Theorem ~\ref{T:#1}}
\newcommand{\refl}[1]{Lemma ~\ref{L:#1}}
\newcommand{\refp}[1]{Proposition ~\ref{P:#1}}
\newcommand{\refc}[1]{Corollary ~\ref{C:#1}}
\newcommand{\refd}[1]{Definition ~\ref{D:#1}}

\newcommand{\refe}[1]{\eqref{E:#1}}

\newenvironment{thm}[1]%
	{ \begin{Thm} \label{T:#1}  \ifShowLabels \TeXref{T:#1} \fi }%
	{ \end{Thm} }

\renewcommand{\th}[1]{\begin{thm}{#1} \sl }
\renewcommand{\eth}{\end{thm} }

\newenvironment{lemma}[1]%
	{ \begin{Lem} \label{L:#1}  \ifShowLabels \TeXref{L:#1} \fi }%
	{ \end{Lem} }
\newcommand{\lem}[1]{\begin{lemma}{#1} \sl}
\newcommand{\elem}{\end{lemma}}

\newenvironment{propos}[1]%
	{ \begin{Prop} \label{P:#1}  \ifShowLabels \TeXref{P:#1} \fi }%
	{ \end{Prop} }
\newcommand{\prop}[1]{\begin{propos}{#1}\sl }
\newcommand{\eprop}{\end{propos}}

\newenvironment{corol}[1]%
	{ \begin{Cor} \label{C:#1}  \ifShowLabels \TeXref{C:#1} \fi }%
	{ \end{Cor} }
\newcommand{\cor}[1]{\begin{corol}{#1} \sl }
\newcommand{\ecor}{\end{corol}}

\newenvironment{defeni}[1]%
	{ \begin{Def} \label{D:#1}  \ifShowLabels \TeXref{D:#1} \fi }%
	{ \end{Def} }
\newcommand{\defe}[1]{\begin{defeni}{#1} \sl }
\newcommand{\edefe}{\end{defeni}}

\newenvironment{remark}[1]%
	{ \begin{Rem} \label{R:#1}  \ifShowLabels \TeXref{R:#1} \fi }%
	{ \end{Rem} }
\newcommand{\rem}[1]{\begin{remark}{#1}}
\newcommand{\erem}{\end{remark}}

\newcommand{\eq}[1]%
	{ \ifShowLabels \TeXref{E:#1} \fi 
	   \begin{equation} \label{E:#1} }
\newcommand{\eeq}{ \end{equation} }

\newcommand{\prf}{ \begin{proof} }
\newcommand{\epr}{ \end{proof} }

\newcommand\alp{\alpha}		
		
\newcommand\gam{\gamma}		\newcommand\Gam{\Gamma}
\newcommand\del{\delta}		\newcommand\Del{\Delta}
\newcommand\eps{\varepsilon}

\newcommand\kap{\kappa}
\newcommand\lam{\lambda}

\newcommand\ome{\omega}		\newcommand\Ome{\Omega}

\newcommand\calA{{\mathcal{A}}}
\newcommand\calB{{\mathcal{B}}}
\newcommand\calC{{\mathcal{C}}}
\newcommand\calD{{\mathcal{D}}}

\newcommand\calF{{\mathcal{F}}}

\newcommand\calH{{\mathcal{H}}}

\newcommand\calK{{\mathcal{K}}}
\newcommand\calL{{\mathcal{L}}}

\newcommand\calO{{\mathcal{O}}}
\newcommand\calP{{\mathcal{P}}}

\newcommand\calR{{\mathcal{R}}}
\newcommand\calS{{\mathcal{S}}}

		\newcommand\bfB{{\mathbf B}}
		
\newcommand\bfd{{\mathbf d}}		
		\newcommand\bfE{{\mathbf E}}
		
		\newcommand\bfG{{\mathbf G}}

		\newcommand\bfP{{\mathbf P}}

\newcommand\bft{{\mathbf t}}		\newcommand\bfT{{\mathbf T}}
		\newcommand\bfU{{\mathbf U}}
		\newcommand\bfV{{\mathbf V}}
		
		\newcommand\bfX{{\mathbf X}}
		\newcommand\bfY{{\mathbf Y}}

\newcommand\RR{\mathbb{R}}

\newcommand\PP{\mathbb{P}}
\renewcommand\AA{\mathbb{A}}

\newcommand\GG{\mathbb{G}}

\newcommand\ZZ{\mathbb{Z}}

\newcommand\CC{\mathbb{C}}

\newcommand\NN{\mathbb{N}}

	\newcommand\grt{{\mathfrak{t}}}

\newcommand\sdp{\times \hskip -0.3em {\raise 0.3ex
\hbox{$\scriptscriptstyle |$}}} 

\newcommand\Hom{\operatorname {Hom}}

\newcommand\id{\operatorname{id}}

\newcommand\Ind{\operatorname{Ind}}

\newcommand\Perv{\operatorname{\Perv}}

\newcommand\SPAN{\operatorname{span}}

\newcommand\Spec{\operatorname{Spec}}

\newcommand\supp{\operatorname{supp}}

\newcommand\ow{{\overline{w}}}

\newcommand\oX{{\overline{X}}}

\newcommand\tilH{{\widetilde{H}}}

\newcommand\tilS{{\widetilde{S}}}

\newcommand\tilW{{\widetilde{W}}}

\newcommand\x{\times}
\newcommand\ten{\otimes}

\renewcommand{\>}{\rangle}
\newcommand{\<}{\langle}

\newcommand{\ra}{\rangle}
\newcommand{\la}{\langle}

\newcommand\sx{\calS(X)}
\newcommand\sdx{\calS'(X)}
\newcommand\sz{\calS^0(X)}
\newcommand\Sc{\calS_c(X)}
\newcommand\vol{\text{vol}}
\newcommand\gv{{\bfG}^{\vee}}
\newcommand\tv{{\bfT}^{\vee}}

\newcommand\yt{Y_{\bfT}}
\newcommand\idt{\bfT(\AA_F)/\bfT(F)}
\newcommand\tw{\grt^w_\RR}
\newcommand\te{\grt^e_\RR}
\newcommand\Eis{\text{Eis}}

\newcommand\alc{{\bf Alc}}
\newcommand\tilPi{{\widetilde \Pi}}

\newcommand\dv{\Del^{\vee}}

\newcommand\av{\alp^{\vee}}

\begin{document}

\title[Schwartz space of the basic affine space]
{On the Schwartz space of the basic affine space}
\author[A.~Braverman and D.~Kazhdan]{Alexander Braverman and David
Kazhdan}
\thanks{The work of both authors was partially supported by the 
Natural Science Foundation}
\address{2-175, Department of Mathematics\\
         Massachusetts Institute of Technology
         77 Massachusetts Ave., Cambridge MA, USA}
\address{Department of Mathematics
        Harvard University
        1 Oxford st. Cambridge MA, USA}

\email{braval@math.mit.edu, kazhdan@math.harvard.edu}

\begin{abstract}
Let $G$ be the group of points of a split reductive 
algebraic  group $\bfG$ over a local field k and let $X=G/U$ 
where $U$ is the group of $k$-points of a maximal unipotent subgroup of 
$\bfG$. In this
paper we construct certain canonical $G$-invariant space 
$\calS(X)$ (called the Schwartz space of $X$) of functions on $X$, 
which is an extension of the space of smooth compactly 
supported functions on $X$. We show that the space of
all elements of $\calS(X)$, which are invariant under the Iwahori subgroup
of $G$ coincides with space generated by the elements of
the so called periodic Lusztig's basis, introduced recently
by G.~Lusztig (cf. \cite{Lu80} and \cite{Lu97}). 
We also give an interpretation of this space in
terms of certain equivariant K-group (this was also done by
G.Lusztig -- cf. \cite{Lu98}). 
Finally we present a global analogue of $\calS(X)$, which
allows us to give a somewhat untraditional treatment of the
theory of principal Eisenstein series.
\end{abstract}
\maketitle

\sec{int}{Introduction}

\ssec{}{}Let $k$ be a non-arcimedian local field and let $\bfG$ be
a split reductive connected algebraic group over $k$. We will assume
also that the derived group $[\bfG,\bfG]$ is simply connected. Let
$\bfX=\bfG/\bfU$ where $\bfU$ is a maximal unipotent subgroup of
$\bfG$ defined over $k$. The group $\bfG$ acts on $\bfX$ on the left. 

In this paper we are going to study
the space of functions on $X\overset{def}=\bfX(k)$ (the set
of $k$-points of $\bfX$). Usually one studies the space $\Sc$
of locally constant, compactly supported functions on $X$. The group
$G=\bfG(k)$ operates naturally on this space. $\Sc$ is usually
called the space of {\it principle series representations}. 
One may observe, however, that for certain purposes $\Sc$ is
not the best possible space of functions on $X$. One way to see
this is the following. 

Let $\bfT$ denote the Cartan group of
$\bfG$. Then $\bfT$ acts on $\bfX$ and this action commutes with
the $\bfG$-action. Hence $T=\bfT(k)$ acts on $\Sc$ and
this action commutes with the $G$-action on $\Sc$. Let $\chi$ be
a character of $T$. Then we denote by $I_\chi$ the space of 
coinvariants of $T$ on $\Sc$ with respect to the character 
$\chi$. Then $I_\chi$ is an admissible representation of
$G$, which is irreducible for generic $\chi$ 
(``$I$'' here stands for ``induced'', since $I_\chi$ is an induced
representation of $G$). Then it is well-known that
there is a natural action $\chi\mapsto w_\bullet \chi$ of $W$ on 
the set of all characters of $T$ such that $I_\chi$ is isomorphic
$I_{w_{\bullet}\chi}$ for generic 
$\chi$ (cf. \refe{dotted} for the definition of
this action). 
However, this is not true for all characters $\chi$. The main
purpose of this paper is ``to correct'' this difficulty. Namely
we define certain bigger space $\sx$ of functions
on $X$ (which we call the {\it Schwartz space} of $X$), which
is $G\x T$-invariant and such
that 

1) $\calS_\chi=I_\chi$ for generic $\chi$, where $\calS_\chi$
denotes the space of $(T,\chi)$-coinvariants in $\sx$.

2) There is a natural action of $W$ on $\sx$, commuting with the
$G$-action, which gives rise to the isomorphisms 
$\calS_\chi{\widetilde \to}\calS_{w_\bullet \chi}$ for any
$w\in W$ and any character $\chi$ of $T$.

\ssec{}{Example: $\bfG=SL(2)$}Suppose now that $\bfG=SL(2)$. In this
case $\bfX=\AA^2\backslash \{ 0\}$ (here $\AA^2$ denotes the affine
plane), hence $X=k^2\backslash \{0\}$. In this case we define
$\sx$ to be space of locally constant, compactly supported
functions on $k^2$. One can check directly that it satisfies
conditions
1 and 2 above. 
\ssec{}{The general case}Unlike the case of $\bfG=SL(2)$, for general
group $\bfG$ the space $\sx$ will {\it not} have any natural
realization
as the space $C^{\infty}_c(\bfY(k))$ of locally constant compactly
supported functions on $\bfY(k)$ for an algebraic variety $\bfY$ over
$k$. So, to define it we need to adopt a rather different strategy.
It was shown by Gelfand and Graev (cf. \cite{Ka} and references 
therein) that the space $L^2(X)$ of $L^2$-functions on $X$ admits 
a natural action of $W$ by  unitary operators 
$\Phi_w$. The operators $\Phi_w$ can be thought of as
generalized Fourier transforms (cf. \refss{fouriertransforms}
for the definition of $\Phi_w$). We define $\sx$ to be equal tthe
subspace
of $L^2(X)$, equal to 
the sum of $\Phi_w(\Sc)$ over all elements $w\in W$ (the sum
is taken inside $L^2(X)$). Despite the fact that this definition
might seem rather artificial, it turns out the the space $\sx$
has a lot of nice properties (the reader should compare our definition
of $\sx$ with the definition of affine Hecke algebras given
in \cite{GKV}).

It would be very interesting to generalize our construction to the
space of functions on $\bfG/\bfU_{\bfP}(k)$ for any
parabolic subgroup $\bfP$ of $\bfG$ (here $\bfU_{\bfP}$ denotes
the unipotent radical of $\bfP$). However, we don't know
how to do this, since the definition of the operators,
analogous to $\Phi_w$ seems to be much
more complicated in the parabolic case.

\ssec{}{}This paper is organized as follows.
\refs{prem} contains some basic notations and preliminaries about the
basic affine space $\bfX$. In \refs{schwartz} we define the Schwartz
space $\sx$ of functions on $X$ and describe explicitely the space
$K$-ivariant vectors there, where $K$ is a maximal compact subgroup
in $G\x T$. In \refs{kl} we study the space of $I$-invariant vectors in
$\sx$, where $I\subset G$ is an Iwahori subgroup. In particular we show
that this space coincides (after suitable identifications) with the space,
spanned by the elements of the so called ``periodic Lusztig's
basis''. In \refs{ktheory} we continue working with the space
$\sx^I$ give another construction of it, using certain equivariant
K-group (thus reproving a result by G.~Lusztig --
cf. \cite{Lu98}). 
In \refs{dualspace} we give certain convenient explicit 
description of the space $\calS'(X)$, which is the dual space to $\sx$
(this description is, in some sense, analogous to the main result in
\cite{GKV}).
Finally in \refs{global} we discuss an analogue of the space
$\sx$ for a global field $F$. In particular we give a somewhat
untraditional treatment of the theory of (principle) Eisenstein
series, which is morally similar to the definition of geometric
Eisentein series (cf. \cite{BG}, \cite{Ga}).

\ssec{}{Acknowledgements}We are very grateful J.~Bernstein for
a lot of useful discussions on the subject and to G.~Lusztig for
numerous helpful conversations and remarks and for letting us know the contents
of \cite{Lu98}. We would also like to thank V.~Ginzburg for pointing
out a similarity between \reft{dualspace} and the main result of
\cite{GKV}. The first author is also indebted to W.~Soergel who
first introduced him to the contents of \cite{Lu97}.

\sec{prem}{Notations and preliminaries}
\ssec{notations}{Notations}Let $k$ be a local non-arcimedean field,
$\calO\subset k$ -- its ring of integers. 
We will denote algebraic varieties over $k$ by bold letters and their
sets of $k$-points by the corresponding ordinary letters.
Let $q$ denote the number of elements in the residue field
of $k$. Let also $\pi$ be a uniformiser of $k$. We will 
suppose that the norm on $k$ is normalized in
such a way that $||\pi||=q^{-1}$.

Let $\bfG$ be reductive, connected, simply connected algebraic
group over $k$. 
We will suppose also that $\bfG$ is split over $k$ and that
its derived group $[\bfG,\bfG]$ is simply connected. 
Let $\bfU$ be a $k$-rational maximal unipotent subgroup of
$k$. Set $\bfX=\bfG/\bfU$. We will refer to $\bfX$ as the 
basic affine space of $\bfG$. The group $\bfG$ acts naturally
on $\bfX$ (on the left). 

Let $\bfB$ be the Borel subgroup of $\bfG$ which contains
$\bfU$ and let $\bfT=\bfB/\bfU$ be the corresponding
Cartan group. Then $\bfT$ also acts naturally on $\bfX$
(on the right) and this action commutes with the $\bfG$-action
on $\bfX$. 

The variety $\bfX$ is quasi-affine. Let ${\overline \bfX}$ denote
its affine closure. One has a natural $\bfG$-equivariant
disjoint decomposition
\eq{}
{\overline \bfX}=\bigcup \bfX_{\bfP}
\end{equation}
where the union is taken over all conjugacy classes of
parabolic subgroups in $\bfG$ and $\bfX_{\bfP}=\bfG/[\bfP,\bfP]$.
Set now 
\eq{}
{\overline \bfX}^{\text{reg}}=\bfX\cup (\bigcup\limits_
{\bfP\ \text{subminimal}}\bfX_{\bfP})
\end{equation}
(by subminimal parabolic we mean a parabolic subgroup $\bfP$
whose Levi group has semisimple rank 1). The variety 
${\overline \bfX}^{\text{reg}}$ is the maximal 
non-singular open subset of ${\overline \bfX}$.

Set $G=\bfG(k),T=\bfT(k),X=\bfX(k)$. 
Let $\Gam=T/\bfT(\calO)$ be the coroot lattice of $\bfG$. 
We will identify every coroot $\alp^{\vee}:~\GG_m\to \bfT$
with its image in $\Gam$ (the latter is well-defined,
since the group $k^{\x}/\calO^{\x}$ is canonically
identified with $\ZZ$). We will denote by
$\Del^{\vee}$ the set of coroots of $\bfG$ and by
$\Del^{\vee}_+$ the set of positive coroots.
We will also denote by $\Gam^+$ the semigroup of positive
elements in $\Gam$. 
We also let $val:~T\to \Gam$ denote the natural projection.

Let $\Gam^{\vee}$ denote the
dual lattice to $\Gam$. This is the lattice of weights of
$\bfG$. Also let ${\widetilde \Gam}$ denote the coweight lattice
of $\bfG$ (thus $\Gam\subset {\widetilde \Gam}$ is a subgroup
of finite index). Let $\Pi\subset \Del_+\subset \Del\subset \Gam^{\vee}$ 
be respectively
the set of simple roots, the set of positive roots and the
set of all roots of $\bfG$. Let also $\rho\in \Gam^{\vee}$ be
the half-sum of all positive roots. 

Let $K=\bfG(\calO)\x \bfT(\calO)$. Then $K$-orbits on $X$ 
are in natural one-to-one correspondence with elements of 
$\Gam$ (in fact any $\bfG(\calO)$-orbit on $X$ is automatically
a $K$-orbit, so in what follows one may use $\bfG(\calO)$ 
instead of $K$). This correspondence can be described in the following
way. Let $\gam\in \Gam$ and let $t_\gam$ be a representative
of $\gam$ in $T$. Then we set $X_\gam=\bfG(\calO)t_\gam\bfT(\calO)/U$.
It is easy to see that
\eq{}
X=\bigcup\limits_{\gam\in \Gam} X_{\gam}
\end{equation}
and that this union is disjoint. 

The variety $\bfX$ admits unique up to a constant top-degree $\bfG$-invariant
differential form, which by \cite{Weil} produces a $G$-invariant
measure on $X$.
It is easy to see that for such a measure one
has 
\eq{}
\vol(X_\gam)=q^{-2\<\gam,\rho\>}\vol(X_0)
\end{equation}
 Therefore, we see that the sum $\sum_{\gam\in \Gam^+}\vol(X_\gam)$
is convergent. We will normalize our measure by requiring
that
\eq{}
\vol(\bigcup\limits_{\gam\in \Gam^+}X_\gam)=1
\end{equation}

Let $C^{\infty}(T)$ denote the space of locally constant functions
on $T$. For any such function $a$ and any $w\in W$ we set
\eq{dotted}
w_{\bullet}(a)(t)=q^{\< val(t), w(\rho)-\rho\>}a(w^{-1}(t))
\end{equation}

This formula defines an action of $w$ on $C^{\infty}(T)$.


\ssec{fouriertransforms}
{Fourier transforms (normalized intertwining operators)}
In what follows we fix a non-trivial character
$\psi:F\to \CC^*$.

Let $L^2(X)$ denote the space of $L^2$-functions on $X$ with
respect to the above measure. By \cite{Ka} this space admits
a natural action of the Weyl group $W$ of $\bfG$ by means
of certain unitary operators $\Phi_w$.
Let $\alp$ be a simple root of $\bfG$ and let $s_\alp\in W$ be
the  corresponding simple reflection. Set $\Phi_\alp=\Phi_{s_\alp}$.
Then the operator $\Phi_\alp$ can be explicitly described as follows.

For a simple root $\alpha$ let $\bfP_\alpha\subset \bfG$ be 
the minimal parabolic
of type $\alpha$ containing $B$. Let $\bfB_\alpha$ be the commutator
subgroup of $P_\alpha$,
and denote $\bfX_\alpha:=\bfG/\bfB_\alpha$. We have an obvious projection
of homogeneous spaces $\pi _\alpha:
\bfX\rightarrow \bfX_\alpha$.
It is a fibration with the fiber $\bfB_\alpha /\bfU=\AA^2-\{0\}$.

Let $\overline \pi _\alpha :{\overline \bfX^\alpha} \to \bfX_\alpha$ 
be the relative
affine completion of the morphism $\pi _\alpha$. (So  $\overline \pi _\alpha$
is the affine morphism corresponding to the sheaf of algebras $\pi_{\alpha *}
(\calO_{\bfX})$ on $\bfX_\alpha$.) Then  $\overline \pi _\alpha$ 
has the structure
of a 2-dimensional vector bundle; $\bfX$ is identified with the complement to
the zero-section in $\overline \bfX^\alpha$. The $\bfG$-action on $\bfX$ 
obviously
extends to $\overline \bfX^\alpha$; moreover, it is easy to see that the
determinant of the vector bundle $\overline \pi _\alpha$ admits a canonical
(up to a constant) $\bfG$-invariant trivialization, i.e.
$\overline \pi _\alpha$ admits unique up to a constant $\bfG$-invariant
fiberwise symplectic form $\ome_\alp$. 
We will fix such a form for every $\alp$.

Obviously $L^2(X)=L^2({\overline X^\alpha})$. Thus we define 
$\Phi_\alp$ to be equal to the Fourier transform in the fibers
of $\overline \pi _\alpha$, corresponding to the identification
of ${\overline X^\alpha}$ with the dual bundle 
by means of $\ome_\alp$.

\sec{schwartz}{The Schwartz space}
\ssec{}{Definition of the Schwartz space}
Let $\Sc$ denote the space of smooth functions with compact support
on $X$ (by a smooth function on $X$ we mean a function which
is invariant under some open compact subgroup of $G\x T$).
Clearly $\Sc\subset L^2(X)$.
Recall that the Weyl group
$W$ acts on $L^2(X)$ by means of the operators $\Phi_w$. 
These operators do not preserve the subspace $\Sc$. We define now
\eq{}
\sx=\sum\limits_{w\in W}\Phi_w(\Sc);\qquad
\sz=\bigcap\limits_{w\in W}\Phi_w(\Sc)
\end{equation}

Both the sum and the intersection are taken in $L^2(X)$.
We have obvious inclusions
\eq{}
\sz\subset \Sc\subset\sx\subset L^2(X)
\end{equation}

\ssec{}{Example}Let $\bfG=SL(2)$. In this case
$X$ is isomorphic to $k^2\backslash \{0\}$ (here $\{0\}\in k^2$
denotes the origin) with the natural action of
$G=SL(2,k)$ and $T=k^{\x}$ on it. We claim now that
$\sx$ in this case coincides with the space 
$C^{\infty}_c(k^2)$ of locally constant compactly supported functions
on $k^2$. 

Indeed, we have the obvious inclusion $\Sc\subset C^{\infty}_c(k^2)$.
On the other hand, $C^{\infty}_c(k^2)$ is invariant under
the operator $\Phi_\alp$ (where $\alp$ is the unique
simple root of $SL(2)$), since in this
case $\Phi_\alp$ just coincides with the symplectic
Fourier transform. Hence $\sx$ is contained in $C^{\infty}_c(k^2)$.

Let us prove the opposite inclusion. 
Let $c_0$ denote the characteristic function of
$\calO^2\subset k^2$. Then it is enough to show that
$c_0\in \sx$ since any $f\in  C^{\infty}_c(k^2)$ can
be decomposed as
\eq{}
f= (f-f(0)c_0)+f(0)\del_0
\end{equation}
and $f-f(0)c_0\in \Sc$.
 
For any $n\in \ZZ$ define functions $c_n\in C^{\infty}_c(k^2)$
and $\del_n\in \Sc$ by
\eq{}
c_n(x)=
\begin{cases}
{q^n\quad \text{if $x\in \pi^n \calO^2$}}\\
{0\quad \text{otherwise}}
\end{cases}
\del_n(x)=
\begin{cases}
{q^n\quad \text{if $x\in \pi^n \calO^2\backslash \pi^{n+1}\calO^2$}}\\
{0\quad \text{otherwise}}
\end{cases}
\end{equation}

Then $\Phi_\alp(c_n)=c_{-n}$. Hence
\eq{}
\begin{align*}
\Phi_{\alp}(\del_n&)=\Phi_{\alp}(c_n-q^{-1}c_{n+1})=\\
&c_{-n}-q^{-1}c_{-n-1}=(1-q^{-2})(c_{-n}-q^{-1}\del_{-n-1})
\end{align*}
\end{equation}
Therefore $c_0=(1-q^{-2})^{-1}(\Phi_\alp(\del_0)+q^{-1}\del_{-1})$
which shows that $c_0\in \sx$.

For general $\bfG$ we cannot give such a simple 
description of $\sx$. However, the following result holds. 
\lem{preliminary}
\begin{enumerate}
\item Any function $f\in \sx$ is well-defined and locally constant
on $\oX^{\text{reg}}$.
\item For any $f\in \sx$ let $\supp_\Gam (f)$ denote the
set of all $\gam\in \Gam$ such that $f|_{X_\gam}\neq 0$. Then there 
exists an element $\gam\in \Gam$ such that
$\supp_\Gam(f)\subset \gam +\Gam^+$.
\end{enumerate}
\elem

\prf First of all, we claim that 2 implies 1. Indeed, we must show that
for any $w\in W$ and any $f\in \Sc$ the function $\Phi_w(f)$ satisfies
1. Let us prove this by induction on $\ell(w)$. So, assume that 
1 holds for some $w\in W$. Take any $\alp\in \Pi$. We must show
that $\Phi_{s_\alp w}(f)=\Phi_\alp(\Phi_w(f))$ satisfies 1. It is easy 
to see from the definition of $\Phi_\alp$ 
that for any simple root $\alp$ of $\bfG$ and any locally constant
function $g$ on $X$, which satisfies both 1 and 2 the function $\Phi_\alp(g)$
is well-defined and satisfies 1. Taking now $g=\Phi_w$ we see that 2
implies 1. 

So, let us prove 2. 
This statement can be reformulated as  follows:
for every $f\in \Sc$ and every $w\in W$ there exists
a locally constant function $a$ on $T$ such that
$\Phi_w(f)=f'*a$ where $f'\in \Sc$ and such that
$a(t)=0$ if $val(t)\not\in \Gam^+$. Let us prove
this by induction on $\ell(w)$. 

First of all, if $\ell(w)=1$, i.e. $w=s_\alp$ is a simple reflection,
then the statement is clear. Moreover, in this case one can
choose $a$ in such a way that $a(t)\neq 0$ only if
$val(t)=n\alp^{\vee}$ for $n\in \NN$. 

Assume now by induction, that we know our statement for some
$w\in W$ and let $\alp\in \Pi$ be a simple root, such that
$\ell(ws_\alp)=\ell(w)+1$. Let $f\in \Sc$. Let $a_1, a_2\in
C^{\infty}(T)$
be two functions, such that $\Phi_\alp(f)=f'*a_1$ and
$\Phi_w(f')=f''*a_2$ such that all the above properties 
hold, i.e. $f',f''\in \Sc$, $a_1(t)\neq 0$   only if
$val(t)=n\alp^{\vee}$ with $n>0$ and $a_2(t)=0$ if $val(t)\not\in \Gam^+$. 
Then $\Phi_{ws_\alp}(f)=f''*(a_2*w_{\bullet}(a_1))$
where $w_{\bullet}(a_1)(t)=q^{\< val(w(t)),\rho\>}a_1(w(t))$.
So, it is enough to check that $a_2*w_{\bullet}(a_1)(t)=0$, if
$val(t)\not\in\Gam^+$. By the definition
\eq{}
a_2*w_{\bullet}(a_1)(t)=\int\limits_{\tau\in T} q^{\< val(w(\tau)),\rho\>}
a_2(t\tau^{-1})a_1(w^{-1}(\tau))d\tau
\end{equation}
We claim now that if $val(t)\not\in \Gam^+$ 
then $a_2(t\tau^{-1})a_1(w^{-1}(\tau))=0$
for any $\tau\in T$. Indeed, if $val(w^{-1}(\tau))\neq n\av$ for some 
$n\in \NN$,
then $a_1(w^{-1}(\tau))=0$. So, suppose that
$val(w^{-1}(\tau))=n\alp^{\vee}$ with $n\in \NN$, i.e.
$val(\tau)=nw(\alp^{\vee})$. But the
condition $\ell(ws_\alp)=\ell(w)+1$ implies that $w(\av)\in \Gam^+$.
Hence $val(\tau)>0$. Therefore $val(t\tau^{-1})\not\in \Gam^+$ which by
our assumptions implies that $a_2(t\tau^{-1})=0$.
This finishes the proof.
\epr
 
\ssec{}{The algebra $\calA$}In what follows we fix a Haar measure on
$T$. We willnormalize it by requiring that the volume of
$\bfT(\calO)$ is equal to 1. Let $\calA$ denote the space of locally
compact, compactly supported functions on $T$. Then $\calA$ has
a natural structure of a commutative algebra with respect to
convolution.

Let $\chi$ be a character of $\bfT(\calO)$. Denote by $\calA_\chi$ the
subspace of $\calA$, consisting of all functions $f\in \calA$ such
that $f(tx)=\chi(t^{-1})f(x)$ for any $t\in \bfT(\calO), x\in T$. 
Then $\calA_\chi$ is a subalgebra of $\calA$, which is in fact
isomorphic to $\CC[\Gam]=\calA_0$ for any $\chi$. Moreover, one
has a direct sum decomposition $\calA=\oplus \calA_\chi$.
Let also $Z_\chi=\Spec \calA_\chi$.

The algebra $\calA$ acts naturally on any smooth representation of 
$T$ and, in particular, on the spaces 
$\sz, \Sc$ and $\sx$. If $V$ is any smooth representation of
$T$ we set 
\eq{}
V^{\chi}=\{v\in V|\ t(v)=\chi(t)v\ \text{for any $t\in \bfT(\calO)$}\}
\end{equation}
The space $V^{\chi}$ has a natural structure of an $\calA_\chi$-module.

Let $\pi:~V_1\to V_2$ be a morphism of smooth $T$-modules.

\defe{} We say that $\pi$ is an isomorphism at the generic
point of $\calA$ if for every character $\chi$ of $\bfT(\calO)$ the
induced morphism $\pi_\chi:~V_1^{\chi}\to V_2^{\chi}$ of
$\calA_\chi$-modules is an isomorphism at the generic
point of $Z_\chi$.
\edefe

\lem{generic}The embeddings $\sz\hookrightarrow \Sc\hookrightarrow \sx$
are isomorphisms at the generic point of $\calA$.
\elem
\prf First of all, it is easy to see that the fact that
the map $\sz\hookrightarrow \Sc$ is an isomorphism at the
generic point of $\calA$ implies that the map
$\Sc\hookrightarrow \sx$ is an isomorphism at the generic point of
$\calA$. 

So, let us show that the map $\sz\hookrightarrow \Sc$ is an
isomorphism
at the generic point of $\calA$. 
Let $\rho$ be a character of $T$. Then $\rho$ defines a point
in $Z_\chi$ where $\chi=\rho|_{\bfT(\calO)}$. 
Let $I_\rho=\Ind_B^G(\rho)$ denote
the
corresponding induced representation (here we regard $\rho$ as
a character of $B$ by means of the identification 
$B/U=T$). Let also
\eq{}
\begin{align*}
I^0_\rho=& \sz/\SPAN\{ f\in\sz| f=t^*(g)-\rho(t)^{-1}g\ \\
&\text{for some }t\in T, g\in \sz\}
\end{align*}
\end{equation}
Then we have the natural map $p_\rho:~I^0_\rho\to I_\rho$. 

Now in order to show the the embedding 
$\sz\hookrightarrow \Sc$ is an isomorphism at the
generic point of $\calA$ it is enough to show that
$p_\rho$ is surjective for generic $\rho$. However,
it is well known that for generic $\rho$ the
$G$-module $I_\rho$ is irreducible. Hence we just
have to show that for generic $\rho$ the morphism
$p_\rho$ is non-zero.

Recall that a character $\Psi$ of $U$ is called 
{\it non-degenerate} if for any $\alp\in \Pi$, the restriction of
$\Psi$ to the one-parametric subgroup $U_\alp$, corresponding to
$\alp$, is non-trivial.
Let now $\Psi$ be a non-degenerate character of $U$.
Set $W_\Psi=\Ind_U^G(\Psi)$ to be the corresponding induced
representation. Then it is well-known that 
for generic $\rho$ one has $\Hom_G(I_\rho,W_\Psi)\neq 0$.
On the other hand by  
\cite{Ka} we have
\eq{}
\Hom_G(\Sc,W_\Psi)\simeq\Hom_G(\sz,W_\Psi)
\end{equation}
where the isomorphism is induced by the embedding $\sz\hookrightarrow
\Sc$.
This implies that $p_\rho\neq 0$ for generic $\rho$.
\epr

\ssec{pairing}{Some pairings}Let $A$ be a commutative algebra and
let $M$ be a module over $A$. Then $M^{\vee}=\Hom_A(M,A)$ also has a natural
structure of an $A$-module. Moreover we have a natural map
$M\to (M^{\vee})^{\vee}$. Recall that an $A$-module 
$M$ is called {\it reflexive} if the above map is an isomorphism.
Clearly any finitely generated free $A$-module is reflexive.
 
Let $V$ be any representation of $G\x T$. We will say that it
is {\it $\calA$-reflexive} if for any open compact
subgroup $H$ of $G$ and any character $\chi$ of $\bfT(\calO)$
the space $V^{H,\chi}$ is reflexive as an
$\calA_\chi$-module (here $V^{H,\chi}$ denotes the space
of all elements of $V$, which are invariant with respect to 
$H$ and on which $\bfT(\calO)$ acts by means of $\chi$). 

For any smooth $G\x T$-module $V$ we define 
\eq{}
V^{\vee}=\sum _{\chi} \bigcup\limits_{H} 
\Hom_{\calA_{\chi}}(V^{H,\chi}, \calA_{\chi})
\end{equation}

Note that $V^{\vee}$ also have a structure of a $G\x T$-module.
It is easy to see that if we have an isomorphism 
$V\simeq V_1^{\vee}$ for some $G\x T$-module $V_1$, which is
isomorphic
at the genric point of $\calA$ to an $\calA$-reflexive module, then
$V$ is also $\calA$-reflexive. 

\th{duality}The $G\x T$-modules $\sz$ and $\Sc$ are
$\calA$-reflexive. Moreover, one has natural identifications
\eq{duality}
\sx^{\vee}\simeq \sz\quad \text{and}\quad
\Sc^{\vee}\simeq \Sc
\end{equation}

These isomorphisms are compatible with the embeddings
$\sz\subset \Sc$ and $\Sc\subset \sx$. Moreover, the first
isomorphism in \refe{duality} is $W$-equivariant (recall, 
that $W$ acts both on $\sz$ and on $\sx$ by $\Phi_w$).
\eth

\prf The fact that $\Sc$ is $\calA$-reflexive follows easily from
the fact that $T$ acts freely on $X$ and from the fact that
for any compact open subgroup $H$ of $G$ the set of $H\x T$-orbits
on $X$ is finite.

Define now a pairing $\<\cdot,\cdot\>:~\Sc\ten \Sc\to \calA$
by
\eq{pairing}
\<f,g\> (t)= q^{-\<val(t),\rho\>}(t^*f, g)_{L^2}
\end{equation}
Here $(\cdot,\cdot)_{L^2}$ denotes the $L^2$ scalar product and
$t^*f(x)=f(tx)$. 

It is easy to see that this pairing is $\calA$-equivariant and that
it defines an isomorphism between $\Sc$ and $\Sc^{\vee}$. 

Note now
that the expression in the right hand side of \refe{pairing} makes
sense for any two smooth functions $f$ and $g$ if we suppose that
$g\in \Sc$. However, if we don't suppose that $f\in \Sc$, 
it might happen that the resulting
function $\la f,g\ra$ does not have compact support and, therefore,
does not lie in $\calA$. 

\lem{moredual}
\begin{enumerate}
\item Suppose that we are given $f\in \sx, g\in \sz$. Then the 
function $\<f,g\>$ defined by \refe{pairing} lies in $\calA$ and
therefore \refe{pairing} defines an $\calA$-equivariant
pairing $\<\cdot,\cdot\>:~\sx\ten \sz\to \calA$.

\item The above pairing is $W$-equivariant.

\item The above pairing gives rise to an isomorphism
$\sz\simeq \sx^{\vee}$.

\item The space $\sz$ is $\calA$-reflexive.
\end{enumerate}
\elem

\prf To prove 1 we must show that  $\<f,g\>$ has compact support.
For this it enough to show that the set of all $\gam\in \Gam$
such that $\<f,g\>|_{\gam \bfT(\calO)}\neq 0$ is finite.
Let us denote this set by $Z_{f,g}$. Then it follows from
\refl{preliminary}(2) that there exists $\gam\in \Gam$
such that $Z_{f,g}$ is contained in $\gam+\Gam^+$. On the
other hand, the same is true also for the set
$Z_{\Phi_w(f),\Phi_w(g)}$. But $Z_{\Phi_w(f),\Phi_w(g)}=w(Z_{f,g})$,
which implies that $Z_{f,g}$ is finite.

The second statement in \refl{moredual} is obvious. Let us 
prove the third one. First of all let us show that our
pairing defines an isomorphism $\sz\simeq \sx^{\vee}$.
Indeed, we have a map $\sz\to \sx^{\vee}$. It is easy to see that
this map is injective. On the other hand, it also easy
to see that we have a natural embedding
$\sx^{\vee}\to \Sc$ (this follows for example from \refl{generic}).
So, let $g\in \sx^{\vee}\subset \Sc$. Then, since $\sx$ is invariant
under the operators $\Phi_w$, it follows that $\Phi_w(g)\subset
\sx^{\vee}$
for any $w\in W$. Hence $\Phi_w(g)\in \Sc$ for any $w\in W$ which
implies that $g\in \sz$. Hence $\sz=\sx^{\vee}$ which implies also
that $\sz$ is $\calA$-reflexive.

\epr 
\epr
\ssec{}{A generalization}Let $Z\subset W$ be any subset. Then we can
define spaces $\calS^Z(X)$ and $\calS^Z_0(X)$ by
\eq{}
\calS^Z(X)=\sum\limits_{w\in Z}\Phi_w(\Sc); \qquad
\calS^Z_0(X)=\bigcap\limits_{w\in Z}\Phi_w(\Sc)
\end{equation}
Then \reft{duality} can be generalized in the following
way (we will need this generalization in \refs{dualspace}).
\th{generalization}The $G\x T$-modules $\calS_0^Z(X)$ is
$\calA$-reflexive. Moreover, one has a natural identification
\eq{duality}
\calS_0^Z(X)\simeq \calS^Z(X)^{\vee}
\end{equation}
\eth

The proof of this result is the same as the proof of \reft{duality}
and it is left to the reader.
\ssec{K-invariant}{Description of $K$-invariant vectors}
Set now $K=\bfG(\calO)\times \bfT(\calO)$. Thus $K$ is a maximal
compact subgroup of $G\times T$. In this subsection
we will give an explicit description of $K$-invariant
elements of $\sx$ (in fact, as it is observed in \refs{prem},
one has $\sx^K=\sx^{\bfG(\calO)}$. So, in what follows, one may
replace $K$ by $\bfG(\calO)$).

\sssec{}{The functions $\del_{\gam}$}
Recall that $\Gam$ denotes the coroot lattice of $\bfG$
and that one has a disjoint decomposition

\eq{}
X=\bigcup\limits_{\gam\in \Gam} X_{\gam}
\end{equation}

Define now a function $\del_{\gam}\in \Sc$ by
$$
\del_{\gam}(x)=
\begin{cases}
{q^{\<\gam,\rho\>}\quad \text{if $x\in X_{\gam}$}}\\
{0\quad \text{otherwise}}
\end{cases}
$$

It is clear that the functions $\del_{\gam}$ form a basis
in $\Sc^K$.

\sssec{}{The $q$-analog of Kostant's function}
Let $\Gam^+\subset \Gam$ denote the semigroup of positive elements
in $\Gam$. Following Lusztig we define a function
$\calK:~ \Gam^+\to \CC$ by 
\eq{}
\calK(\gam)=\sum\limits_{P\in \calP_\gam}q^{-|P|}
\end{equation}
where $\calP_\gam$ is the set of all representations
of $\gam$ as a sum of positive coroots and if
$P\in \calP_\gam$ then we denote by $|P|$ the number of
summands in $P$.

The function $\calK(\gam)$ may be viewed as a $q$-analog of
Kostant's partition function.

\sssec{}{The functions $c_{\mu}$} Let us now define a function
$c_{\mu}$ on $X$ by putting
\eq{emu}
c_{\mu}=\sum\limits_{\gam\in \Gam^+}\calK(\gam) \del_{\mu+\gam}
\end{equation}

\noindent
{\it Remark.}\quad The reader should compare our definition of
the function $c_{\mu}$ with the computation of stalks of intersection
cohomology sheaves on Drinfeld spaces, given in \cite{FFKM}. 
\th{K-invariant}

\begin{enumerate}
\item $c_{\mu}\in \sx$ for any $\mu\in \Gam$
\item The functions $c_{\mu}$ form a basis in $\sx^K$
\item $\Phi_w(c_\mu)=c_{w(\mu)}$ for any $w\in W, \mu\in \Gam$
\end{enumerate}
\eth

\prf Let us first prove 3. First of all we claim that 
$e_\mu\in L^2(X)$ for any $\mu\in \Gam$.
Indeed, one has
\eq{norm}
||c_\mu||^2_{L^2}=\vol(X_0)^2 
\sum\limits_{\gam \in \Gam^+}\calK(\gam)^2
\end{equation}
So, we have to show that the sum in the right hand side of \refe{norm}
converges.
But it is easy to see that
\eq{}
\sum\limits_{\gam\in \Gam^+} \calK(\gam)=\frac{1}{ (1-q^{-1})^
{|\Del_+^{\vee}|}}
\end{equation}

Hence $\sum\limits_{\gam\in \Gam^+}\calK(\gam)$ is convergent and therefore
$\sum\limits_{\gam \in \Gam^+}\calK(\gam)^2$ is also convergent.

Thus the expression $\Phi_w(c_\mu)$ is well defined. 
Clearly in order to prove \reft{K-invariant} one may
assume without loss of generality that $w=s_\alp$ for
some $\alp\in \Pi$ and $\mu=0$.

Let us introduce some notations.
Let $\gam\in \Gam$. 
Then we can define two operators $\tilH_\gam$ and $H_{\gam}$
on $\bfT(\calO)$-invariant functions on $X$: $\tilH_\gam$
is just the operator of shift by $\gam$ and 
$H_\gam=q^{\<\gam,\rho\>}\tilH_\gam$. Note that
$\tilH_{\gam_1}\tilH_{\gam_2}=\tilH_{\gam_1+\gam_2}$ and 
$H_{\gam_1}H_{\gam_2}=H_{\gam_1+\gam_2}$.

The following lemma is proved by a direct calculation (which reduces
essentially to $G=SL(2)$.
\lem{}
\eq{salp}
\Phi_\alp
(\del_0)=\frac{1-q^{-1}H_{-\alp^{\vee}}}{1-q^{-1}H_{\alp^{\vee}}}
(\del_0)
\end{equation}

Here $\frac{1}{1-q^{-1}H_{\alp^{\vee}}}$ is understood as a formal
sum 
\eq{}
\frac{1}{1-q^{-1}H_{\alp^{\vee}}}=\sum q^{-n}H_{n\alp^{\vee}}
\end{equation}
\elem

On the other hand, we claim that the following
equality holds
\eq{obraschenie}
(\prod\limits_{\alp^{\vee}\in \Del^{\vee}_+}q-H_{\alp^{\vee}})(c_\mu)=
q^{-|\Del^{\vee}_+|}\del_\mu
\end{equation}
where $H_\alp$ denotes the operator of shift by $\alp$. To see that
\refe{obraschenie} holds it is enough to note the following.
Let ${\widehat \CC[\Gam]}$ denote the completion of $\CC[\Gam]$
consisting of all infinite sums $\sum a_\gam H_\gam$ satisfying
the following condition:

$\bullet$ for every $\mu\in \Gam$ the set
\eq{}
\{ \gam\in \Gam|\ a_\gam\neq 0\ \text{and } \gam\not\in \mu +\Gam^+\}
\end{equation}
is finite.

It is obvious now that in ${\widehat \CC[\Gam]}$ we have
the following equality:
\eq{bre}
\frac{1}{\prod\limits_{\alp^{\vee}\in \Del^{\vee}_+}q-H_{\alp^{\vee}}}=
q^{|\Del_+^{\vee}|}\sum_{\gam\in \Gam^+}\calK(\gam)H_\gam
\end{equation}

On the other hand \refe{bre} implies \refe{obraschenie}
 because
of the fact that $H_\gam(\del_\mu)=\del_{\mu+\gam}$.
But \refe{salp} and \refe{obraschenie} imply together that
$\Phi_\alp(c_0)=c_0$ which is precisely what we need to prove.

Let us prove 1 and 2 now. Let $N$ denote the linear span
of the functions $c_\mu$. It is clear that $N$ is invariant
under the action of $\Gam$ (by shifts). Also \reft{K-invariant}(3)
implies
that $N$ is also invariant under $W$. 

We claim now that for any
$\mu\in \Gam$ one has $\del_\mu\in N$. 
Indeed, after one expands \refe{obraschenie},
one gets an expression of $\del_\mu$ as a finite sum of
$c_{\mu'}$'s with some coefficients, since 
$H_\gam(c_\mu)=c_{\mu+\gam}$.

Let $N_c=\Sc^K$. It follows from the above that $N_c\subset N$.
Hence, in order to show that $N=\sx^K$ (which is precisely
what we have to prove), it is enough
to show that the spaces $\Phi_w(N_c)$ generate
$N$.

Let $\bfT^{\vee}$ denote the dual torus of $\bfT$. By definition, this
is an algebraic torus over $\CC$, whose character group is
$\Gam$. The group $W$ acts naturally on $\bfT^{\vee}$. Let 
$\CC(\bfT^{\vee})$ denote the space of regular functions on $\bfT^{\vee}$.
It can be naturally identified with the group ring of $\CC[\Gam]$
of $\Gam$. 

It follows from \reft{K-invariant}(3) that the assignment
$c_\mu\mapsto \mu$ extends to an isomorphism
$N\simeq \CC(\bfT^{\vee})$ of $\Gam\rtimes  W$-modules, where we let
$\Gam$ act on $N$ by means of the operators $H_\gam$. 
Moreover, it follows from \refe{obraschenie} that
under this isomorphism the space $N_c$ gets identified
with the ideal of $\CC(\bfT^{\vee})$, generated by
the function $\bfd_q=\prod_{\alp\in \Del_+^{\vee}}(q-\alp^{\vee})$
(here we consider $\alp^{\vee}$ as a function of $\bfT^{\vee}$ with
values in $\GG_m$).

Let $\bfV\subset \bfT^{\vee}$ denote hypersurface in
$\bfT^{\vee}$ defined by the equation $\bfd_q=0$.
\lem{empty} One has 
\eq{empty}
\bigcap\limits_{w\in W}w(\bfV)=\emptyset
\end{equation}
\elem

\prf By Hilbert Nullstellensatz it is enough to prove that for any
$\bft^{\vee}\in \bfT^{\vee}$ there exists $w\in W$ such that
$f(w(\bft^{\vee}))\neq 0$. However, for every $\bft^{\vee}\in \tv$
there exists $w\in W$ such that
\eq{}
|\gam(w(\bft^{\vee}))|\leq 1
\end{equation}
for every $\gam\in \Gam^+$
(this means that $w(\bft^{\vee})$ lies in the ``antidominant
chamber''), which proves the lemma.
\epr
On the other hand, it follows immediately from \refe{empty}
that the sum all $w$-translates of the ideal, generated by
$\bfd_q$ in $\CC(\bfT^{\vee})$ is equal to $\CC(\bfT^{\vee})$, 
which finishes the proof.
\epr
\sec{kl}{Relation with generic Kazhdan-Lusztig polynomials}
In \refs{kl} and \refs{ktheory} we suppose that $\bfG$ is
semisimple.
\ssec{}{The affine Hecke algebra}
Let $\tilW=\Gam\rtimes W$. This is a Coxeter group. We will denote by
$\tilPi$ the corresponding set of generators.
 By definition this is an algebra over the
ring $\CC[v,v^{-1}]$ of Laurent polynomials in an indeterminate
$v$ with generators $T_w$ for $w\in \tilW$. We will denote 
$T_{s_\alp}$ simply by $T_\alp$. 
 These generators satisfy
the following relations:
\eq{}
\begin{align}
&(T_\alp+v^{-1})(T_\alp-v)=0\quad \\
&T_{ww'}=T_wT_{w'}\quad \text{if 
$\ell(ww')=\ell(w)+\ell(w')$}
\end{align}
\end{equation}

(here $\ell:\tilW\to \ZZ^+$ is the length function).
We will write $\calH_q$ for the specialization of
$\calH$ at $v=q^{\frac{1}{2}}$. Let also $\calH^f\subset \calH$ (resp.
$\calH^f_q\subset \calH_q$) denote the finite Hecke algebra, i.e.
the subalgebra of $\calH$ spanned by the $T_w$ for $w\in W$.

Let $I\subset G$ denote the Iwahori subgroup.
It is well-known that $\calH_q$ is isomorphic to the algebra of 
smooth compactly supported functions on $G$ which are 
$I$-bi-equivariant. This isomorphism can be characterized
as follows. It is well known that the set $I\backslash G/I$ is
isomorphic to $\tilW$. For any $w\in \tilW$ let $\chi_w$ denote the
characterictic function of the corresponding double coset.
Then under the above isomorphism the element
$(-1)^{\ell(w)}T^{-1}_{w^{-1}}$ 
goes to
$q^{-\frac{1}{2}\ell(w)}\chi_w$.

Therefore, it follows that  
if $V$ is any smooth representation of $G$ then the algebra
$\calH_q$ acts naturally on $V^I$.

\ssec{}{$\Sc^{I\x \bfT(\calO)}$ and the periodic Hecke module}
Let $\alc$ denote the collection of all alcoves for the
group $\tilW$ in the space $\grt_{\RR}$ -- the 
real Lie algebra of $\bfT$. The set $\alc$ admits two
commuting (resp. left and right actions) of the group
$\tilW$. Let also $\calC^+\in \grt_\RR$
denote the dominant Weyl chamber. For any $s\in \tilS$ we will
denote by $H_s$ the corresponding affine hyperplane in
$\grt_{\RR}$. For any $s\in \tilS$ we denote by
$H_s^+$ (resp. $H_s^-$ the set of all points $x\in \grt_\RR$ which
are on the same side of $H_s$ as the set $\gam+\calC^+$ (resp.
$\gam-\calC^+$ for some $\gam\in \Gam$. We also denote by
$d:\alc\to \ZZ$ the corresponding length function
(cf. \cite{Lu80}(Section 1.4)). 

The set $\alc$ has two natural (respectively left and right) 
actions of $\tilW$ which commute. We will denote them by
$A\mapsto wA$ (resp. $A\mapsto Aw$).  

Following Lusztig we define for any 
$A\in \alc$ a subset $\calL(A)\subset \tilPi$ in the following
way. We say that $\alp\in \calL(A)$ if and only if the two conditions
below are satisfied:

$\bullet$ $A\subset E_\alp^+$.

$\bullet$ $As_\alp\subset E_\alp^-$ (here $As$ denotes the right
action
of $s_\alp$ on $A$).
   
In \cite{Lu80} and \cite{Lu97} G.~Lusztig considered the
{\it periodic Hecke module} $M_c$ over $\calH$. 
By definition this module is 
spanned by elements $\{A\in \alc\}$ over the ring $\CC[v,v^{-1}]$
and the action of 
$H$ is defined in the following way:
\eq{taction}
T_{\alp}A=
\begin{cases}
s_\alp A\quad \text{if $s\not\in \calL(A)$}\\
s_\alp A+(v-v^{-1})A\quad \text{if $s\in \calL(A)$}
\end{cases}
\end{equation}
(thus, if we specialize $v$ to $1$ we obtain the natural
$\CC[\tilW]$-action on $\SPAN(A|\ A\in \alc)$). 

The module $M_c$ also admits a natural action of the group $\Gam$,
which commutes with the $\calH$-action. This action is defined
by $\gam(A)=A+\gam$.

Let also $M_{\geq}$ denote a completion of $M_c$, consisting
of all (possibly infinite) sums $\sum m_A A$ such that
the following property holds:

$\bullet$ there exists $\gam\in \Gam$, such that
$m_A\neq 0$ implies that $A\subset \gam+\calC^+$.

In \cite{Lu97} G.~Lusztig has defined 
for any $\alp\in \Pi$ certain $\calH$-linear operator
$\theta_\alp: M_c\to M_{\geq}$ in the following way. 
Let $A\in \alc$. Define a sequence $A^n, n\geq 0$ by
the consitions that $A^n$ lie in the same ``$\alp$-strip''
as $A$ and

$\bullet$ $A^0=As_\alp$

$\bullet$ $d_\alp(A^n,A^{n+1})=1$, where $d_\alp$ denotes the
corresponding ``partial''
distance function (cf. \cite{Lu97}).
Then 
\eq{}
\theta_\alp(A)=v^{-1}A^0+\sum\limits_{n=1}^{\infty}(v^{-n+1}-v^{-n-1})A^n
\end{equation}
It clear, that when we specialize $v$ to 1 the operator
$\theta_{\alp}$ reduces just to an obvious operator on
$M$, which comes from the right action of $W$ on $\alc$.

The set $\alc$ is also in natural one-to-one correspondence
with the set of $I\x\bfT(\calO)$-orbits on $X$ (in fact, as before,
any $I$-invariant set is automatically $\bfT(\calO)$-invariant).
So, we may write 
\eq{}
X=\bigcup\limits_{A\in \alc} X_A
\end{equation}
where $X_A$ is the corresponding orbit. For any
$A\in \alc$ we define a function $\del_A\in \Sc^{I\x \bfT(\calO)}$
by
\eq{}
\del_A(x)=
\begin{cases}
(-q^{1/2})^{d(A)}\quad\text{if $x\in X_A$}\\
0\qquad\text{otherwise}
\end{cases}
\end{equation}

The following lemma can be proved by a direct computation and it is
left to the reader.
\lem{mc}
Let $M_{c,q}$ denote the specialization of $M_c$ at $v=q^{-\frac{1}{2}}$.
Then the assignment $A\mapsto \del_A$ extends to an isomorphism
$\rho:M_{c,q}{\widetilde \to} \Sc^{I\x \bfT(\calO)}$
of $\calH_q\ten \CC[\Gam]$-modules. Moreover,
$\rho\circ\theta_\alp\circ\rho^{-1}=\Phi_\alp$ for any $\alp\in \Pi$.
\elem

\ssec{}{}In \cite{Lu97}(Theorem 12.2) G.~Lusztig has constructed
certain remarkable elements $A^\sharp\in M_{\geq}$. 
These elements satisfy the following conditions: first of all, 
$A^\sharp=A+v^{-1}\sum_B \ZZ[v^{-1}] B$. Also $A^\sharp$ is 
``self-dual'' in a certain sense (cf. \cite{Lu97}).

Let $M_d=\SPAN(A^{\sharp}|\ A\in \alc)$. Then it follows from
Theorem 12.2 in \cite{Lu97} that the operators $\theta_\alp$ are
well-defined on $M_d$ and together define an action of the 
group $W$ on $M_d$. We will denote the corresponding automorphisms
of $M_d$ by $\theta_w$.

\prop{lusztig}One has
\eq{lusztig}
\sum\limits_{w\in W}\theta_w(M_c)=M_d
\end{equation}
\eprop

\prf Let 
\eq{}
\bfd_v=\prod\limits_{\alp^{\vee}\in \Del^{\vee}_+} (v^2-\alp^{\vee})\in 
\calR_{T\x \CC^*}
\end{equation}

Then it follows from \cite{Kato}(Section 3) (cf. also Conjecture
13.11 in \cite{Lu97} and the remark afterwards) that
$\bfd_v M_d\subset M_c$. Hence \refe{lusztig} follows from the
fact the ideal generated by the elements $\{w(\bfd_v)\}_{w\in W}$
in $\calR_{\bfT^{\vee}\x \CC^*}$ contains 1, which is proved exactly
in the same way as \refl{empty}
\epr

\cor{lusztig}Let $M_{d,q}$ denote the specialization of $M_d$ at
$v=q^{\frac{1}{2}}$. Then there is an isomorphism of 
$\calH_q\ten \CC[\Gam\rtimes W]$-modules 
\eq{}
M_{d,q}=\sx^{I\x \bfT(\calO)}
\end{equation}
\ecor

This follows immediately from \refp{lusztig} and \refl{mc}.

\sec{ktheory}{Relation with equivariant $K$-theory.}
\ssec{}{}The results of this section were also obtained  earlier
by G.~Lusztig, using slightly different methods (cf. \cite{Lu98}).
\ssec{}{}Let $\gv$ the 
Langlands dual group to $\bfG$ and let $\tv\subset \gv$ be
its Cartan subgroup. Let also $\calB$ denote the flag variety of
$\gv$. One has a natural identification 
$\Gam\simeq \Hom(\tv,\CC^*)$. 

Let $\calK=K_0^{\tv\x \CC^*}(\calB)$ denote the equivariant
$K$-group of $\calB$ with respect to the $\tv\x \CC^*$ (where
$\CC^*$ acts trivially on $\calB$). This is a module over
the ring $\CC[v,v^{-1}]$ which can be identified with the ring
of characters of $\CC^*$. We will write $\calK_q$ for the
specialization
of $\calK$ at $v=q$. 

By \cite{KL} the algebra
$\calH$ acts naturally on $\calK$ (as in \cite{Lu97} we twist the
action defined in \cite{KL} by means of the automorphism
$\bullet:\calH\to\calH$ sending $T_w$ to
$(-1)^{\ell(w)}T^{-1}_{w^{-1}}$). 

On the other hand, the space $\calK$ admits also a natural action
of the group $\Gam\rtimes W$, which commutes with the action
of $\calH$.  This action may be described as follows. 

First
of all, since $\Gam$ is identified with the lattice of
characters of $\tv$, it follows that $\Gam\simeq K_0^{\tv}(pt)$, and
therefore, $\Gam$ acts on $K_0^{\tv}(Y)$ for any $\tv$-variety
$Y$. 

Suppose now that $Y$ is actually a $\gv$-variety.
Then the Weyl group $W$ acts naturally on $K_0^{\tv}(Y)$ in 
the following way. Let $N(\tv)$ denote the normalizer of
$\tv$ in $\gv$. Then $W$ can be identified with $N(\tv)/\tv$.
Let $w\in W$ and let $\ow$ be a representative of $w$ in
$N(\tv)$. Let also $\calF$ be $\tv$-equivariant coherent sheaf
on $Y$. Then $(\ow^{-1})^*(\calF)$ also acquires a natural structure
of a $\tv$-equivariant coherent sheaf, which is does not depend
on the choice of $\ow$ up to an isomorphism. Hence
$(\ow^{-1})^*(\calF)$
defines an element in $K_0^{\tv}(Y)$. It is easy to see that
this element depends only on the image of $\calF$ in
$(\ow^{-1})^*(\calF)$. Hence we get an action of $W$
on $(\ow^{-1})^*(\calF)$.

Applying this now to $Y=\calB$ we see that the
space $K_0^{\tv\x\CC^*}(\calB)$ acquires an action of the
algebra $\calH\ten \CC[\Gam\rtimes W]$, where
$\CC[\Gam\rtimes W]$ denotes the group algebra of
$\Gam\rtimes W$.

\ssec{}{}By \cite{Kato} (cf. also \cite{Lu97}) $M_c$ is isomorphic
to $\calH$ as a $\calH\ten \calR_{\tv\x \CC^*}$-module. Under this
isomorphism the element $1\in \calH$ goes to $A^+\in M_c$, where
$A^+$ is the unique alcove, which is contained in $\calC^+$ and
contains $0$ in its closure.

Let 
\eq{}
2\rho^{\vee}=\sum\limits_{\alp^{\vee}\in\dv_+}\alp^{\vee}
\end{equation}
Consider now the set of fixed points of $\tv$ on $\calB$. These
fixed points are in natural one-to-one correspondence with
the elements of $W$. For any $w\in W$ we will denote
by $y_w\in\calB $ the corresponding fixed point. Let $\kap_w$ denote
the skyscraper sheaf at $y_w$. We let $\tv$ act on the fiber 
of $\kap_w$ by the character $w(2\rho^{\vee})$. We will denote
the image of $\kap_w$ in $\calK$ by the same symbol.
Set $\kap=\kap_e$.

Following
\cite{KL} we define a homomorphism $\zeta:M_c\to \calK$ of
$\calH\ten\calR_{\tv\x \CC^*}$-modules by sending $A^+$ to $\kap$.
By \cite{KL}(Theorem 6.2) $\zeta$ is well-defined and it becomes
an isomorphism after tensoring with the ring of fractions
of $\calR_{\tv\x \CC^*}$. More precisely the following result
holds:
\th{kl}(cf. Theorem 6.2 in \cite{KL}) There exists $N\in \NN$ such that
$\bfd_v^N \calK\subset \zeta(M_c)$ (in fact, in this it is known that
one can take $N=1$ but we are not going to use this).
\eth

Hence we see that the image of $\zeta^{-1}$ (which {\it a priori})
maps
$\calK$ into some localization of $M_c$) lies in $M_{\geq}$.

The main result of this section is the following 
\th{k-theory} 
\begin{enumerate}
\item The morphism $\zeta$ is ``$W$-equivariant'', i.e. 
\eq{phizeta}
\theta_\alp=\zeta^{-1}\circ(s_\alp)^*\circ\zeta
\end{equation}
\item 
\eq{}
\sum\limits_{w\in W}w^*(\zeta(M_c))=\calK
\end{equation}
\end{enumerate}
\eth

\cor{klusztig}The map $\zeta$ extends to an isomorphism between $M_d$
and $\calK$.
\ecor

This follows immediately from \reft{k-theory} and \refp{lusztig}.

\cor{ks}$\calK_q$ is naturally isomorphic to
$\calS(X)^{I\x\bfT(\calO)}$
as a $\calH\ten \CC[\tilW]$-module.
\ecor

Note that we do not need \refc{lusztig} in order to prove \refc{ks}.

\ssec{}{Proof of \reft{k-theory}}

As in \cite{Lu97} we can regard both $M_c$ and $\calK$ as sections
of some vector bundles $\bfV\to \tv\x \CC^*$ and $\bfV'\to \tv\x\CC^*$.
The map $\zeta$ induces a meromorphic isomorphism between $\bfV$ and
$\bfV'$. Now both the left hand side and the right hand side
of \refe{phizeta} are meromorphic maps from $\bfV$ to itself, which
induce the same map $s_\alp$ on the base. 
On the other hand, it is well known that the fiber of $\bfV$ at
the generic point of $\tv\x \CC^*$ is an irreducible representation
of $\calH$. Hence the left hand side and the right hand side of
\refe{phizeta} differ by multiplication by certain meromorphic
function $f$ on $\tv\x \CC^*$.

On the other hand we know that 
$\theta_\alp^2= (\zeta^{-1}\circ s_\alp^*\circ\zeta)^2=\id$.
This implies that 
\eq{f=1}
f(\bft,\lam)f(s_\alp(\bft),\lam)=1
\end{equation}

We claim now that 
$f(\bft,1)=1$ for every $\bft\in \tv$ (it is easy
to see that the divisor $v=1$ is not a singular divisor for $f$, hence
the restriction of $f$ to $v=1$ makes sense).
To prove this we must show that $\zeta$ commutes with the
action of $W$ when we specialize to $v=1$.
We will denote the corresponding specializations by
$\zeta^1$, $M_c^1$, $\calK^1$ etc. (Note that for
$v=1$ the operators $\theta^1_\alp$ preserve $M_c^1$ and
therefore define an action of $W$ there).
By \refe{taction} one has $T^1_w A^{+}=A^+ w$. It is easy to see that
one also has $T^1_\alp\kap=s_\alp(\kap)$.  Indeed, fix $\alp\in \Pi$. 
Let $\pi_\alp:~\calB\to \calB_\alp$ denote the projection
from $\calB$ to the corresponding partial flag variety. Let
also $\Ome_\alp$ denote the corresponding sheaf of relative
differential forms. 
Then
\eq{}
T^1_\alp(\calF)=\pi^*_\alp\circ(\pi_\alp)_*(\calF)-
\pi^*_\alp\circ(\pi_\alp)_*(\calF\ten \Ome_\alp)-\calF
\end{equation}
Let
$\calP=\pi_\alp^{-1}(\pi_\alp(y_e))=\pi_\alp^{-1}(\pi_\alp(y_{s_\alp})$.
Then $\calP$ is $\tv$-stable and therefore the structure
sheaf $\calO_{\calP}$ acquires the natural $\tv$-equivariant
structure. Also $\calP$ is $\tv$-equivariantly isomorphic
to the projective line $\PP^1$, where we let $\tv$ act on $\PP^1$
by means of the character $\alp^{\vee}:\tv\to \CC^*$, composed
with the standard $\CC^*$-action on $\PP^1$. 
For any $\gam\in {\widetilde\Gam}$ and $n\in \ZZ$, such that
 $n+\<\alp,\gam\>$ is even, we will denote by
$\calO(n)_{\calP}^\gam$ the image in $\calK^1$ of the
bundle $\calO(n)_{\calP}\ten \CC_\gam$ (where $\CC_\gam$ is the
one-dimensional representation of a covering of $\tv$, on which the 
Lie algebra $\grt^{\vee}$ acts by
means of $d\gam$). Thus
$$
T^1_\alp(\kap)=\calO_{\calP}^{2\rho^{\vee}}+
\calO_{\calP}^{2\rho^{\vee}-2\av}-\kap$$.
But it is easy to see that (note also that 
$2\rho^{\vee}-2\av=s_\alp(\rho^{\vee})$)
\eq{}
\calO_{\calP}^{2\rho^{\vee}}-\kap=\calO_{\calP}^{s_\alp(2\rho^{\vee})}-
\kap(y_{s_\alp})=\calO(-1)_{\calP}^{\rho^{\vee}}
\end{equation}

Hence $T^1_\alp(\kap)=\kap(y_{s_\alp})=s_\alp(\kap)$.
Thus $T^1_w\kap=w(\kap)$ for
every $w\in W$.
This implies that $f(\bft,1)=1$. 

On the other hand, it follows from the above, that
$f$ must be of the form $f=\chi \bfd_v^N$, where $\chi$ is a regular
invertible function on $\tv\x \CC^*$, for some integer $N$.
Thus 
$$
1=fs_\alp^*(f)=\chi s_\alp^*(\chi)\bfd_v^N s_\alp^*(\bfd_v^N)
$$
which implies that $\bfd_v^N s_\alp^*(\bfd_v^N)$ is regular and
invertible.
It is easy to see that this may happen only for $N=0$.

Thus $f=\chi$ is regular and invertible.
But then $f$ is just a scalar multiple of a character
of $\tv\x \CC^*$ and therefore \refe{f=1} and the fact that
 $f(\bft,1)=1$ for every $\bft\in \tv$ implies that $f$ must be
identically equal to 1. Hence \refe{phizeta} holds.

We claim now that \reft{k-theory}(2) follows immediately from 
\reft{kl}. Indeed, arguing in the same way as in \refl{empty}
we can show that the ideal in $\calR_{\tv\x \CC^*}$ generated
by the functions $\{w(\bfd_v^N)\}_{w\in W}$ contains 1, which 
implies \reft{k-theory}(2) in view of \reft{kl}. $\square$


\sec{dualspace}{Description of the dual space}
\ssec{}{}Let $\sdx$ denote the dual space to $\sx$. Let also 
$\calD(X)$ denote the space of distributions on $X$ (i.e. dual space
to $\Sc$). Clearly $W$ acts on $\sdx$ and we will denote the corresponding
operators also by $\Phi_w$. Also one has a surjective
map $r:~\sdx\to \calD(X)$. The main result of this
section is the following theorem, which will be used in the next
subsection
where we discuss
a global analogue of $\sx$.

First of all let us introduce some notations. For every simple
root $\alp$ of $G$ let 
\eq{}
\calS^0_{\alp}(X)=\{ f\in \Sc|\ \text{such that\ } \Phi_{\alp}(f)\in \Sc\}
\end{equation}

\th{dualspace}
The space $\sdx$ is isomorphic to the space of all collections
$\{\lam_w\}_{w\in W}$ where $\lam_w\in \calD(X)$ satisfying the
following condition
\eq{condition}
\lam_{s_\alp w}(f)=\lam_w(\Phi_{s_\alp}(f))\quad \text{for any
$w\in W, \alp\in \Pi, f\in \calS^0_\alp(X)$}
\end{equation}
Moreover, this isomorphism can be described by sending
$\lam\in \sdx$ to $\{\lam_w\}$, where 
$\lam_w=r(\Phi_w(\lam))$. 
\eth
\prf

\rem{} Note the similarity between \reft{dualspace} and the main
result
of \cite{GKV}.
\erem
Choose $\alp\in \Pi$. Let us denote by $W^\alp$ the subset
of $W$, consisting of all minimal representatives of $W/\{1,s_\alp\}$.

Consider the complex $K^2\to K^1\to K^0$ where

$\bullet$ $K^0=\sx$

$\bullet$ $K^1=\oplus_{w\in W}\Sc$

$\bullet$ $K^2=\oplus_{\alp\in \Pi ,w\in W^{\alp}}\calS^0_\alp(X)$

where the differential $d:~K^1\to K^0$ is given by
\eq{}
d(\{f_w\})=\sum\limits_{w\in W}\Phi_w(f_w)
\end{equation} 
and the
differential $d:~K^2\to K^1$ is given by the following formula.
Suppose that we are given an element $k\in K^2$, whose
$(w,\alp)$-component
is equal to $f$ for some $\alp\in \Pi, w\in W^\alp$ and 
$f\in \calS^0_\alp(X)$. Then $d(k)$ is defined by
\eq{}
d(k)_w=f,\ d(k)_{ws_\alp}=-\Phi_\alp(f),\ d(k)_{w'}=0\ \text{otherwise}
\end{equation}
It is easy to see now that the assertion of \reft{dualspace} is
equivalent to the fact that the above complex is exact in the middle
term.
\ssec{complex}{The complex $C^*(M)$} 
Let $\grt_{\RR}=\Gam\ten \RR$ be a real Cartan
algebra of $\bfG$. Define a convex polytope $B$ in $\grt_{\RR}$ in 
the following way. Choose a regular point $a\in \grt$ and
define $B$ to be the convex hull of the set $w(a), w\in W$. 

For any $I\subset \Pi$ let $W^I$ denote the set of all minimal
representatives of $W/W_I$, where $W_I$ is the subgroup of
$W$ generated by $s_i$ for $i\in I$.
Then $B$ admits a natural cell decomposition  such that the cells
of dimension $i$ are parametrized by all pairs $(I\subset
\Pi, w\in W^I)$
such that $|I|=i$. Let $C^*(B)$ denote the ``augmented'' cohomology
complex associated with this cell decomposition. This means
that $C^i(B)$ is the space with a basis parametrized by cell
of dimension $i$ of $B$ and $C^0(B)=\CC$ with the natural
differential $d:~C^i\to C^{i+1}$. The complex 
$C^*(B)$ is obviously acyclic, since $B$ is contractible.

Let $M$ be any $W$-module. Then
we can define a  complex $(C^*(M),d_M)$ in the following manner.
First of all we set $C^i(M)= M\ten C^i(B)$. Let now $m\ten c\in
C^i(M)$.
Suppose that 
\eq{}
d(c)=\sum\limits_{w\in W, \alp \in \Pi\setminus I}
c_{w,I\cup \{\alp\} }
\end{equation}
(here $c_{w,I\cup \{\alp\} }$ denotes the corresponding component
of $d(c)$).
Then we define
\eq{}
d_M(m\ten c)= \sum\limits_{w\in W, \alp\in \Pi\setminus I}
 s_\alp(m)\ten c_{w,I\cup \{\alp\} }
\end{equation}

The complex $(C^*(M), d_M)$ is still acyclic because
its differential is conjugate to $\id\ten d$.

\ssec{}{End of the proof}Recall that we have to prove that
the complex $K^2\to K^1\to K^0$ is acyclic in the middle term.
But it follows from \reft{duality} and \reft{generalization}
that this is equivalent to showing that 
$(K^0)^{\vee}\to (K^1)^{\vee}\to (K^2)^{\vee}$ is acyclic in the
middle term (we use here the fact that $K^1$ and $K^2$ are 
$\calA$-reflexive; note that we don't need to know that
$K^0$ is $\calA$-reflexive).

However, it follows again from  \reft{duality} and
\reft{generalization}
that 
the latter complex embeds naturally into the complex
$C^*(\sx)$, which is acyclic by \refss{complex}.
Let now $x\in (K^1)^{\vee}$ be such that $d(x)=0$.
Then there exists some $y\in \sx$ such that 
$d_{\sx}(y)=x$. However
$d_{\sx}(y)=\oplus_{w\in W}\Phi_w(y)$. On the other hand
it follows once more from \reft{duality} that 
$(K^1)^{\vee}=\oplus_{w\in W}\Sc$. Hence we see that $\Phi_w(y)\in
\Sc$
for every $w\in W$. This means that $y\in \sz=(K^0)^{\vee}$ which
finishes the proof.

\epr
\sec{global}{The global space}
\ssec{}{Remarks on arcimedean places}We now want to describe
some global analogue of the space $\sx$. In order to do that 
we will need a definition of $\sx$ in the case where the
base local field $k$ is arcimedian. So, let
$\bfG$ be a semisimple connected simply connected algebraic
group over $k=\RR,\CC$. We will assume in fact
that $k=\RR$ (which can be done without loss of generality) and
that $\bfG$ is quasi-split over $k$. Clearly the space
$L^2(X)$ still makes sense and one still has an action of 
$W$ on this space defined exactly as before. Let $\calD$ denote the algebra of
all global sections of the sheaf of regular differential
operators on $\bfX$. Thus we define
\eq{}
\sx=\{ f\in L^2(X)|\ \text{such that $d(f)\in L^2(X)$ for every 
$d\in \calD$}\}
\end{equation}
It is easy to see that  $\sx$ is invariant under the operators
$\Phi_w$.

In the sequel we will need the following estimate on the growth
rate of elements of $\calS$. Choose a maximal compact subgroup
$K_\RR$ of $G=\bfG(\RR)$. One can choose a Borel subgroup
$\bfB$ of $\bfG$ in such a way that
$G=K_\RR\cdot B$, where $B=\bfB(\RR)$. Hence the natural
map $\phi:~\bfT(\RR)=T\to X$ becomes surjective after we 
quotient out $X$ by $K_\RR$. Let also 
$\xi_\RR:T\to \grt_\RR$ denote the homomorphism of
abelian groups, defined by the condition
\eq{}
\<\lam^{\vee},\xi_\RR(t)\>=\ln |\lam^{\vee}(t)|
\end{equation}
for every weight $\lam^{\vee}\in \Gam^{\vee}$.

\lem{growth} Let $f\in \sx$. Then for every $n\in \NN$ there
exists a real number $c=c(n)$ such that
\eq{}
|f(\phi(t))|< c |\<\rho,\xi_\RR(t)\>|^{-n}
\end{equation}
\elem

The proof is left to the reader.
\ssec{}{More notations}
Let now $F$ be a global field and let $\bfG$ be again a semisimple
simply connected algebraic group defined over $F$. Let also
$\AA_F$ denote the ring of adeles of $F$. 
Let $v$ be a place of $F$ and let $F_v$ denote the corresponding
local field.

In this case the variety $X$ is defined over $F$. We will suppose that
the symplectic forms $\ome_\alp$ (cf. \refss{notations}) are also
$F$-rational.

Also we will fix a non-trivial character $\psi:\AA_F/F\to \CC^*$.
This character defines a non-trivial character $\psi_v$ of $F_v$
for every place $v$.
All the local Fourier transforms that we are going to consider
will be with respect to the characters $\psi_v$.

\ssec{}{The global space}Let $v$ be a place of $F$. Then we may consider the
space $\calS_v=\calS(\bfX(F_v))$. By \reft{K-invariant}
the space $\calS_v$ has distinguished $K_v=\bfG(\calO_v)\x \bfT(\calO_v)$
invariant vector $c_0^v$ (we do not define $e^v_0$ when 
$v$ is archimedean).

\defe{}The space $\calS$ is the restricted tensor product of the 
spaces $\calS_v$ (over all places $v$ of $F$) with respect to the
vectors $c_0^v$. We will also endow $\calS$ with a topology, which
comes from the discrete topology on $\calS(\bfX(F_v)$ for non-arcimedean
$v$ and from Frechet topology on $\calS(\bfX(F_v)$ for 
archimedean $v$. 
\edefe

Since the vector $c_0^v\in \calS_v$ is $W$-invariant, it follows that
the $W$-action extends (diagonally) to the space $\calS$. 
As before for any $w\in W$ we let $\Phi_w$ denote the corresponding
operator. On the other hand, it is clear that $\calS$ also admits a natural
smooth action of the group $\bfG(\AA_F)\x \bfT(\AA_F)$ compatible
in the obvious sense with the action of $W$.

\ssec{}{The functional $\eps_0$}We claim now that there exists a natural
morphism from $\calS$ to the space of smooth functions on $\bfX(\AA_F)$.
Indeed if $f=\otimes f_v\in \calS$ and $x=\{x_v\}\in \bfX(\AA_F)$
then we define
\eq{defi}
f(x)= \prod f_v(x_v)
\end{equation}
It is easy to see that the product in the left hand side of
\refe{defi} is actually finite (this follows from the fact
that $x_v\in X_{0,v}$ for almost any $v$).

We define now
\eq{}
\eps_0(f)=\sum\limits_{x\in \bfX(F)\subset \bfX(\AA_F)}f(x)
\end{equation}
 
It is clear that the functional $\eps_0$ on $\calS$ is $\bfG(F)\x \bfT(F)$-
invariant, but not $W$-invariant. Our next goal is to produce some natural
$W$-invariant modification of $\eps$.
\ssec{}{The spaces $\calS_{c,v}$}
Fix a place $v$ of $F$. Let $\calS_{c,v}$ denote the subspace of
$\calS$ spanned by all functions of the form
$\otimes f_{v'}$ where $f_v\in \calS_c(\bfX(F_v))$. 
Also, for any $\alp\in \Pi$ let 
\eq{}
\calS^0_{\alp,v}=\{f\in \calS_{c,v}|\ \Phi_{\alp}(f)\in \calS_{c,v}\}
\end{equation}

Set also $\calS^0_\alp=\sum_v \calS^0_{\alp,v}$.

\prop{invariance}
For any simple root $\alp$ and any $f\in \calS^0_{\alp}$ one
has
\eq{}
\eps_0(f)=\eps_0(\Phi_{\alp}(f))
\end{equation}
\eprop
\prf Suppose first that $G=SL(2)$. In this case $\calS$ is the space of
smooth compactly supported functions on $\AA^2_F$ and the space 
$\calS^0_{\alp}$ consists of all functions $f$ such that 
\eq{bred}
f(0)=\Phi(f)(0)=0
\end{equation}
(here we have set $\Phi=\Phi_\alp$).
Also one has
\eq{eps}
\eps_0(f)=\sum\limits_{x\in F^2,x\neq 0}f(x)
\end{equation}

It follows now from the Poisson summation formula that
for any $f\in \calS$ one has
\eq{poisson}
\sum\limits_{x\in F^2}f(x)=\sum\limits_{x\in F^2}\Phi(f)(x)
\end{equation}

Now it is easy to see that \refe{bred}, \refe{eps} and \refe{poisson}
imply together that $\eps_0(f)=\eps_0(\Phi(f))$ for any $f\in \calS^0_\alp$
which proves \refp{invariance} for $\bfG=SL(2)$.

In the general case \refp{invariance} follows from 
the same arguments as above and the following result.

\prop{}
Let $p:\bfE\to \bfY$ be a symplectic
vector bundle, defined over a global field $F$.
 Let $f$ be any 
function on $\bfE(\AA_F)$ which satisfies the 
following condition:

for every $y\in \bfY(\AA_F)$ the restriction of $f$ to $p^{-1}(y)$ lies 
in the space of Schwartz-Bruhat functions of $p^{-1}(y)$.
 
The the Fourier transform  $\Phi(f)$ is well-defined, and, moreover,
one has
\eq{}
\sum\limits_{x\in E(F)}f(x)=\sum\limits_{x\in E(F)} \Phi(f)(x)
\end{equation}
\eprop
\epr

\th{automorphic}
There exists unique $W$-invariant functional $\eps$ on $\calS$ such that
for any place $v$ of $F$ the restriction of $\eps$ to $\calS_{c,v}$ coincides
with the restriction of $\eps_0$ to $\calS_{c,v}$. The functional $\eps$ is
$\bfG(F)\x \bfT(F)$-invariant.
\eth

\prf The uniqueness of $\eps$ is obvious. Namely, let $f\in \calS$. 
Choose a place $v$ of $F$. Then there exist functions
$f_w\in \calS_{c,v}$ such that 
\eq{}
f=\sum\limits_{w\in W}\Phi_w(f_w)
\end{equation}

We set now 
\eq{defeps}
\eps(f)=\sum\limits_{w\in W}\eps(f_w)
\end{equation}
(It is clear that if $\eps$ exists then it must satisfy \refe{defeps}).

In order to show that $\eps$ is well-defined, we must show that
the right hand of \refe{defeps} does not depend on the choice of the
functions $f_w$ and on the choice of $v$. However, the first statement
follows immediately from \reft{dualspace} and the second statement follows
from the following lemma, whose proof is exactly the same as the 
proof of \reft{dualspace}
\prop{} Let $v_1,v_2$ be two places of $F$.
Let $V=\calS(\bfX(F_{v_1}))\ten \calS(\bfX(F_{v_2}))$. The
group $W$ acts on $V$ diagonally.
Set also 
\eq{}
V_c=\calS_c(\bfX(F_{v_1}))\ten \calS(\bfX(F_{v_2}))+
\calS_c(\bfX(F_{v_1}))\ten \calS(\bfX(F_{v_2}))
\end{equation}
Set $V^0_\alp=V_c\cap \Phi_\alp(V_c)$.
Let also $V'_c$ denote the dual space to $V_c$. 

Suppose now that we are given a collection $\{\lam_w\}_{w\in W}$,
where
$\lam_w\in V_c'$, 
such that $\lam_w(f)=\lam_{s_\alp w}(\Phi_\alp(f_{s_\alp w}))$
for any $f\in V^0_\alp$. Then 
$\{\lam_w\}$ gives rise to a well defined functional on $V$.
\eprop
\epr
\ssec{}{}Let $Y_{\bfT}=(\bfG(\AA_F)\x \bfT(\AA_F))/(\bfG(F)\x
\bfT(F))$
and let $C^{\infty}(\yt)$ denote the space of smooth functions on
$\yt$. Thus we may define may $\bfG(\AA_F)\x \bfT(\AA_F)\rtimes W$-
equivariant map $\eta:~\calS\to C^{\infty}(\yt)$ by
\eq{eta}
\eta(f)(g,t)=\eps((g\x t)(f))
\end{equation}

where $(g\x t)(f)(x)=f(g^{-1}x t^{-1})$. 
Let also $\eta_0$ be the corresponding map, defined as in \refe{eta}
but with $\eps$ replaced by $\eps_0$.

\ssec{}{}Let now $\chi$ be a character of $\idt$. We wish to define
a map $E_\chi$ from $\calS$ to the space of functions on
$\bfG(\AA_F)/\bfG(F)$ by putting
\eq{composition}
E_\chi(f)(g)=\int\limits_{\idt} \eta(f)(g,t)\chi(t) dt
\end{equation}

However, the integral in the right hand side of 
\refe{composition} does not have to converge and we will
have to regularize it. We will do that
only for {\it regular} characters $\chi$ (cf. \refd{regular} below).
In fact, for any regular character we will define a morphism
of $\bfG(\AA_F)\x \bfT(\AA_F)$-modules $E_\chi:~\calS\to 
C^{\infty}(\bfG(\AA_F)/\bfG(F)$, where $\bfT(\AA_F)$ acts
on $C^{\infty}(\bfG(\AA_F)/\bfG(F)$ by means of $\chi$, which
``morally'' will be the regularization of the above integral.

Let $\chi$ be a character of $\idt$ and let $w\in W$.
Then we set $w_{\bullet}\chi =||\cdot||^{-1}w(\chi||\cdot||)$.

\defe{regular} We will say that $\chi$ is {\it regular} if 
for any corroot $\alp^{\vee}:~\GG_m\to \bfT$ and any $w\in W$, 
the induced character
$w_{\bullet}\chi\circ \alp^{\vee}$ of $\AA_F^*/F^*$ is non-trivial.
\edefe

Note that the set of regular characters is invariant under the action
of $W$.

\ssec{}{Regularization of the integral for regular characters}
Let $\AA_F^*/F^*$ be the idele class group of $F$. Then we
have the norm map $||\cdot||:~\AA_F^*/F^*\to \RR^+$.

Consider now the group $\idt$. Recall that $\grt_\RR$ denote the
real Cartan algebra of $\bfG$ (i.e. $\grt_\RR=\Gam\ten \RR$. 
Then there exists unique
homorphism of abelian groups $\xi:~\idt\to \grt_\RR$ satisfying
\eq{xi}
\< \lam^{\vee},\xi(t)\>=ln ||\lam^{\vee}(t)||
\end{equation}

for any $t\in \idt$ and $\lam^{\vee}\in \Gam^{\vee}$ (note that 
$\lam^{\vee}$ defines a map $\idt\to \AA_F^*/F^*$ which, abusing
the notations, we denote by the same letter).

Let 
\eq{}
\grt_\RR=\bigcup\limits_{w\in W}\tw
\end{equation} 
be the decomposition of $\grt_\RR$ into Weyl
chambers. 

\prop{convergence}For any $f\in \calS$ the integral
\eq{}
\int\limits_{\xi^{-1}(\te)} \eta_0(f)(g,t) \chi(t)dt
\end{equation} 
is absolutely convergent. Here $e\in W$
denotes the unit element.
\eprop

\prf \refp{convergence} follows immediately from the following
lemma.
\lem{}Let $f\in \calS$.
Then for any $g\in \bfG(\AA_F)/\bfG(F)$ and $n\in \NN$ there
exists a real number $c=c(g,n)$ such that
\eq{}
|\eta_0(f)(g,t)|< c |\<\rho,\xi(t)\>|^{-n}
\end{equation}
\elem

The lemma follows easily from \refl{growth} and \refl{preliminary}.
\epr

Suppose now that $\chi$ is a regular character of $\idt$. 
Then we can define the map $E_\chi:~\calS\to
C^{\infty}(\bfG(\AA_F)/\bfG(F))$
by

\eq{}
E_\chi(f)(g)=
\sum\limits_{w\in W}\int\limits_{\xi^{-1}(\te)}
\eta_0(\Phi_w)(f)(g,t)w_{\bullet}(\chi(t)) dt
\end{equation}

\prop{echi}
\begin{enumerate}
\item $E_\chi$ is a $\bfG(\AA_F)\x W$-equivariant map (where
$W$ acts trivially on $C^{\infty}(\bfG(\AA_F)/\bfG(F))$.

\item Let $\calS_\chi$ denote the space of $(\bfT(\AA_F),\chi)$-
coinvariants on $\calS$, i.e.
\eq{}
\calS_\chi=\calS/{\overline 
\SPAN ( \chi(t)g-t(g)|\ g\in \calS, t\in \bfT(\AA_F))}
\end{equation}
(here  the bar denotes the closure of the corresponding subspace). 
Then $E_\chi$ descends to a well-defined map
$\Eis_\chi:~\calS_\chi\to C^{\infty}(\bfG(\AA_F)/\bfG(F))$.
\end{enumerate}
\eprop

\rem{} We will see some connection between $\Eis_\chi$ and
Eisenstein series in the next subsection.
\erem

\prf 
The second statement of \refp{echi} is clear from the
definitions. Let us prove the first one. 

The $\bfG(\AA_F)$-equivariance of $E_\chi$ is obvious.
Let us prove that $E_\chi$ is $W$-equivariant, i.e. 
that $E_\chi(f)=E_\chi(\Phi_w(f))$ for every $f\in \calS, w\in W$. 
Clearly, it is enough to show this when $w=s_\alp$ is a simple
reflection. Hence \refp{echi} follows from the
following lemma.

\lem{echi-alp} 
For every $f\in \calS$ and any regular character $\chi$ one has
\eq{}
\int\limits_{\xi^{-1}(\te)}
\eta_0(f)(g,t)(\chi(t)) dt=
\int\limits_{\xi^{-1}(\te)}
\eta_0(\Phi_\alp(f))(g,t)s_\alp(\chi(t)) dt
\end{equation}
\elem

\prf (Of \refl{echi-alp}).
Recall that for every simple root $\alp$ we have
defined the $\bfG(\AA_F)\x \bfT(\AA_F)$-submodule 
$\calS^0_\alp$ of $\calS$ which is invariant under $\Phi_\alp$.
Moreover, by \refp{invariance}, one has

\eq{}
\eta_0(f)=\eta_0(\Phi_\alp(f))
\end{equation}
for every $f\in \calS^0_\alp$. Therefore, \refl{echi-alp} holds
for every $f\in \calS^0_\alp$. Hence \refl{echi-alp} follows from the
following easy lemma.
\lem{coinvariants}
The natural map $(\calS^0_\alp)_\chi\to \calS_\chi$ is an isomorphism.
\elem

\epr
\epr
\ssec{}{Connection with usual Eisenstein series} Let 
\eq{}
\pi_\chi=\Ind_{\bfB(\AA_F)}^{\bfG(\AA_F)}\chi
\end{equation}
where we
regard $\chi$ as a character of $\bfB(\AA_F)$ by means
of the identification $\bfB(\AA_F)/\bfU(\AA_F)=\bfT(\AA_F)$. Then by
\cite{Lan} one has the Eisenstein series map 
$\Eis'_\chi:~\pi_\chi\to C^{\infty}(\bfG(\AA_F)/\bfG(F))$, which
is meromorphic in $\chi$ (in particular, it is well-defined for 
generic $\chi$).

Choose now a finite set of places $P$ of $F$ such that
$P$ contains all places, where $\chi$ is ramified and all
arcimedean places of $F$. Then we can define a map
$i_{P,\chi}:~\pi_\chi\to \calS_\chi$ in the following way.

Let $K_P=\prod_{v\not\in P}\bfG(\calO_v)$. The representation
$\pi_\chi$ is generated by its $K_P$-invariant vectors.
Hence it is enough to construct a map
$\pi_\chi^{K_P}\to \calS_\chi^{K_P}$ which commutes with the
action of the corresponding Hecke algebra.
 
By the definition of an induced representation, $\pi$ is
a quotient of the space $C^{\infty}_c(\bfX(\AA_F))$ of smooth
compactly supported functions on $\bfX(\AA_F)$. 
Moreover, $\pi_\chi^{K_P}$ is a quotient of 
$(C^{\infty}_c(\bfX(\AA_F))^{K_P}$. On the other hand, 
we have the natural morphism $l_P:~\calS_\chi^{K_P}\to 
(C^{\infty}_c(\bfX(\AA_F))^{K_P}$ sending $\otimes_v f_v$
to 
\eq{}
\bigotimes\limits_{v\in P}f_v \otimes \bigotimes_{v\not\in P}
(\prod\limits_{\alp^{\vee}\in \Del_+^{\vee}}(1-q_v^{-1}H_{\alp^{\vee}}))(f_v)
\end{equation}
where $q_v$ is the number of elements in the residue field of 
$F_v$ and $H_{\alp^{\vee}}$ is as in \refss{K-invariant}. 
Thus composing $l_P$ with 
$\Eis_\chi$ we get the map
$\Eis_{\chi,P}:~\pi_\chi\to C^{\infty}(\bfG(\AA_F)/\bfG(F))$.

Let now
\eq{}
L_P(\chi)=\prod\limits_{\alp^{\vee}\in\Del_+^{\vee}}
\prod\limits_{v\not\in P}(1-q_v^{-1+\<\alp^{\vee},\rho\>}
{\chi_{\alp^{\vee}}(\pi_v)})^{-1}
\end{equation}

be the corresponding $L$-function, which is a meromorphic function of $\chi$.
Here the notations are the following: $\chi_\alp$ denotes the character
of $\AA_F^*/F^*$ obtained by composing $\chi$ with the homomorphism
$\AA_F^*/F^*\to \idt$, coming from $\alp^{\vee}:~\GG_m\to \bfT$. Also
$\pi_v$ is a uniformiser of $F_v$, viewed as an element of $\AA_F^*$.
 
The following proposition is obvious from the definitions.
\prop{L-function} 
\eq{}
\Eis_{\chi,P}=L_P(\chi)\Eis'_\chi
\end{equation}
\eprop

\rem{} One can view \refp{L-function} as the reason that $L$-function
appears in the functional equation, satisfied by Eisenstein series
(cf. \cite{Lan}). Note that in the case when $F$ is a functional
field and $\chi$ is everywhere unramified one can take $P=\emptyset$.
In this case $\Eis_{\chi,P}$ was constructed in \cite{Laum}, \cite{Ga}
and \cite{BG} by geometric methods. It would be interesting to
understand a direct connection betweeen the construction in
{\it loc. cit.} and the construction presented above.
\erem

\end{document}